\documentclass[12pt]{amsart}
\usepackage{amssymb}
\DeclareMathOperator{\id}{Id}
\DeclareMathOperator{\im}{Im}
\DeclareMathOperator{\image}{image}
\DeclareMathOperator{\crit}{Crit}
\DeclareMathOperator{\car}{Char}
\DeclareMathOperator{\DR}{DR}
\DeclareMathOperator{\supp}{Supp}
\DeclareMathOperator{\gr}{gr}
\newcommand{\Xt}{{\wt X}}
\newcommand{\pDR}{{{}^p\!\DR}}
\newcommand{\DRa}{\DR^{\rm an}}
\newcommand{\pDRa}{\pDR^{\rm an}}
\newcommand{\pCC}{{{}^p\!\CC}}
\newcommand{\phip}{{{}^p\!\phi}}
\newcommand{\pch}{{\check p}}

\newcommand{\qc}{{\check q}}

\newcommand{\ooplus}{\mathop\oplus\limits}
\newcommand{\ootimes}{\mathop\otimes\limits}
\newcommand{\lefpar}{\left(}
\newcommand{\rigpar}{\right)}

\newcommand{\lcr}{\left[\!\left[}
\newcommand{\rcr}{\right]\!\right]}

\let\dpl\displaystyle
\let\wh\widehat
\let\wt\widetilde

\def\qqbox#1{\quad\mbox{#1}\quad}
\let\wwwh\wh
\let\wwh\wh
\def\implique{\Longrightarrow}

\newcommand{\defin}{:=}

\newcommand{\eg}{{\it e.g}}
\newcommand{\cf}{\emph{cf}}

\newcommand{\ie}{{\it i.e}}

\newcommand{\T}{\S\kern .15em }
\newcommand{\ptbl}{.\kern .15em }
\newcommand{\bbullet}{{\scriptscriptstyle\bullet}}

\DeclareMathAlphabet{\mathcalmaigre}{U}{eus}{m}{n}

\def\AA{\mathbf{A}}
\def\CC{\mathbb{C}}
\def\NN{\mathbb{N}}
\def\ZZ{\mathbb{Z}}
\def\PP{\mathbb{P}}
\def\bH{\boldsymbol{H}}
\def\bL{\boldsymbol{L}}
\def\bR{\boldsymbol{R}}

\def\cD{\mathcal{D}}
\def\cE{\mathcal{E}}

\def\cH{\mathcal{H}}

\def\cM{\mathcal{M}}

\def\cO{\mathcal{O}}

\def\cR{\mathcal{R}}

\def\cU{\mathcal{U}}

\def\cX{\mathcalmaigre{X}}

\def\ccF{\mathcalmaigre{F}}

\def\Afuc{\check\AA\!^1}

\newdimen\lengtharrow
 \newbox\exponantbox \newbox\indicebox
\def\dimmax#1#2{\ifdim#1<#2 #2\else #1\fi}

\def\arrowr#1#2%
    {\setbox\exponantbox=\hbox{$#2$}
      \setbox\indicebox=\hbox{$#1$}
       \lengtharrow=\dimmax{\wd\indicebox}{\wd\exponantbox}
        \ifdim \lengtharrow=0pt
															\mathrel{\smash{\mathop{\hbox to 10mm{\rightarrowfill}}}}
         \else \ifdim\lengtharrow<6mm
                \lengtharrow=10mm
                 \else \advance\lengtharrow by 4mm
                  \fi
					\mathrel{\smash{\mathop{\hbox to\lengtharrow{\rightarrowfill}}\limits%
                         ^{\textstyle#2}_{\textstyle#1}}}
									\fi}

\def\arrowl#1#2%
    {\setbox\exponantbox=\hbox{$#2$}
      \setbox\indicebox=\hbox{$#1$}
       \lengtharrow=\dimmax{\wd\indicebox}{\wd\exponantbox}
        \ifdim \lengtharrow=0pt
															\mathrel{\smash{\mathop{\hbox to 10mm{\leftarrowfill}}}}
         \else \ifdim\lengtharrow<6mm
                \lengtharrow=10mm
                 \else \advance\lengtharrow by 4mm
                  \fi
					\mathrel{\smash{\mathop{\hbox to\lengtharrow{\leftarrowfill}}\limits%
                         ^{\textstyle#2}_{\textstyle#1}}}
									\fi}

\def\MRE#1{\arrowr{}{#1}}

\newcommand{\isom}{\stackrel{\sim}{\longrightarrow}}

\newtheorem{theo}{Theorem}
\newtheorem*{cor}{Corollary}
\newtheorem*{prop}{Proposition}
\newtheorem*{lem}{Lemma}

\begin{document}

\title{On a twisted de~Rham complex}
\author{Claude Sabbah}
\address{UMR 7640 du CNRS,
Centre de Math\'ematiques,
\'Ecole polytechnique,
F--91128 Palaiseau cedex,
France}
\email{sabbah@math.polytechnique.fr}
\date{}
\begin{abstract}
We show that, given a projective regular function $f:X\rightarrow \CC$ on a smooth quasi-projective variety, the corresponding cohomology groups of the twisted de~Rham complex $(\Omega_{X}^{\bbullet},d-df\wedge)$ and of the complex $(\Omega_{X}^{\bbullet},df\wedge)$ have the same dimension. We generalize the result to de~Rham complexes with coefficients in a mixed Hodge Module.
\end{abstract}
\maketitle

\section*{Introduction}
\subsection{}
Let $\wt X$ be a smooth projective variety over $\CC$ and $(\Omega_\Xt^{\bbullet},d)$ be the complex of algebraic differential forms. Hodge theory and GAGA theorem of Serre (see also \cite{D-I87} for an algebraic argument, or \cite{Illusie96} for other references) show that the hypercohomology spaces on $\wt X$ of both complexes $(\Omega_\Xt^{\bbullet},d)$ and $(\Omega_\Xt^{\bbullet},0)$ have the same dimension (this follows from the degeneracy at $E_1$ of the spectral sequence $\mbox{Hodge}\implique\mbox{de Rham}$).

\subsection{}
Denote by \label{numfilt}$\cO_\Xt$ the sheaf of regular functions on $\Xt$ and by $\cD_\Xt$ the sheaf of differential operators with coefficients in $\cO_\Xt$. More generally, let $(\wt M,F)$ be a mixed Hodge Module on $\wt X$ as defined by M. Saito \cite[\T4]{MSaito87}, where $F_\bbullet\wt M$ is in particular a good filtration (increasing and exhaustive) of the $\cD_\Xt$-Module $\wt M$. The (algebraic) de~Rham complex $\DR(\wt M)=(\Omega_\Xt^\bbullet\otimes_{\cO_\Xt}\wt M,\nabla)$ is naturally filtered using $F$ (see \cite{Brylinski83}): the degree-$\ell$ term of $F_k\DR(\wt M)$ is $\Omega_{\Xt}^{\ell}\otimes F_{k+\ell}\wt M$, that is also denoted by
\begin{eqnarray*}
F_k\DR(\wt M)&=&(\Omega_{\Xt}^{\bbullet}\otimes_{\cO_\Xt}F_k[-\bbullet]\wt M,\nabla).
\end{eqnarray*}
The associated graded complex is equal to the complex $(\Omega_{\Xt}^{\bbullet}\otimes_{\cO_\Xt}\gr^F\!\wt M,\gr^F\!\nabla)$, where $\gr^F\!\nabla=\oplus_k[\nabla]_{k}^{k+1}$ and $[\nabla]_{k}^{k+1}$ is the degree-$1$ $\cO_\Xt$-linear morphism induced by $\nabla$
\begin{eqnarray*}
\Omega_{\Xt}^{\ell}\otimes_{\cO_\Xt}\gr_k^F\wt M&\longrightarrow & \Omega_{\Xt}^{\ell+1}\otimes_{\cO_\Xt}\gr_{k+1}^F\wt M.
\end{eqnarray*}

The degeneracy at $E_1$ (see \cite[(4.1.3)]{MSaito87}) now implies that the hypercohomology spaces on $\Xt$ of the complexes $\DR(\wt M)$ and $(\Omega_\Xt^\bbullet\otimes_{\cO_\Xt}\gr^F\!\wt M,\gr^F\!\nabla)$ have the same dimension.

\subsection{}\label{numnot}
Let $\wt f:\wt X\rightarrow \PP^1$ be a morphism of algebraic varieties and let $f:X\rightarrow \AA^1$ be its restriction over the affine line $\AA^1$. Thus, $X$ is quasi-projective and $f$ is projective.

\begin{theo}\label{thmglob}
Let $(M,F)$ be a mixed Hodge Module on $X$. The hypercohomology spaces on $X$ of the complexes $(\Omega_X^\bbullet\otimes_{\cO_X} M,\nabla-df\wedge)$ and $(\Omega_X^\bbullet\otimes_{\cO_X}\gr^F\!\!M,\gr^F\!\nabla-df\wedge)$ have the same (finite) dimension.
\end{theo}

\subsection{Remark}\label{remphi}
If \label{numrem}$\phi_f$ denotes the vanishing cycle functor as defined by Deligne \cite{Deligne73} (see also \cite{K-S90}) and $\DRa$ the analytic de~Rham functor, it is well-known (\cf. \T\ref{num21}) that one has
\begin{eqnarray*}
\dim\bH^i\lefpar X,(\Omega_X^{\bbullet}\otimes_{\cO_X} M,\nabla-df\wedge)\rigpar &=&\sum_{c\in\CC}\dim\bH^{i-1}\lefpar f^{-1}(c),\phi_{f-c}\DRa(M)\rigpar.
\end{eqnarray*}

\subsection{}
One \label{numkont}may, after M. Saito \cite[(4.1.2)]{MSaito86}, apply this result to $M=\cO_X$ equipped with its natural filtration. One then has $\gr^F\!\!M=\cO_X$ and $\gr^F\!\nabla=0$. In this way one recovers the result of M. Kontsevitch and S. Barannikov which motivated this work.

\begin{cor}[Kontsevitch and Barannikov]
The hypercohomology spaces on $X$ of the complexes $(\Omega_X^\bbullet,d-df\wedge)$ and $(\Omega_X^\bbullet,df\wedge)$ have the same (finite) dimension.
\end{cor}

\subsection{}
One also \label{numlog}may, after M. Saito \cite[(4.2.2)]{MSaito86}, apply this result to $M=\cO_X[*D]$, if $D$ is a divisor of $X$. Assume that $D$ is a divisor with normal crossings. One gets (see \T\ref{numrem2})

\begin{cor}
The hypercohomology spaces on $X$ of the complexes $(\Omega_X^\bbullet(\log D),d-df\wedge)$ and $(\Omega_X^\bbullet(\log D),df\wedge)$ have the same (finite) dimension.
\end{cor}

\subsection{}
Let \label{numan}now $f:\cX\rightarrow \CC$ be a holomorphic function on a complex analytic manifold. Denote by $\cO_\cX$ the sheaf of holomorphic functions on $\cX$ and by $\cD_\cX$ the corresponding sheaf of differential operators. Let $(\cM,F)$ be a well-filtered coherent $\cD_\cX$-module, $\DRa(\cM)=(\Omega_\cX^\bbullet\otimes_{\cO_\cX}\cM,\nabla)$ its (analytic) de~Rham complex, which is naturally filtered, and $(\Omega_\cX^\bbullet\otimes_{\cO_\cX}\gr^F\!\!\cM,\gr^F\!\nabla)$ the associated graded complex.

It will be convenient to use sometimes the perverse shift convention
$$
\pDRa(\cM)\defin\DRa(\cM)[\dim\cX]=(\Omega_\cX^{\bbullet+\dim\cX}\otimes\cM,\nabla)
$$
and $\phip_f\ccF\defin\phi_f\ccF[-1]$ for a complex $\ccF$ with constructible cohomology on $\cX$.

\begin{theo}\label{thmloc}
Let $(\cM,F)$ be a mixed Hodge Module (according to M. Saito {\rm\cite[\T2.d]{MSaito87}}) such that the restriction of $f$ to the support of $\cM$ is a projective morphism in a neighbourhood of $f^{-1}(0)$. Then for all $i\in\ZZ$ one has
\begin{eqnarray*}
\dim\bH^{i-1}\lefpar f^{-1}(0),\phi_f\DRa(\cM)\rigpar &=&\dim\bH^i\lefpar f^{-1}(0),(\Omega_\cX^{\bbullet} \otimes_{\cO_\cX}\gr^F\!\!\cM,\gr^F\!\nabla-df\wedge)\rigpar.
\end{eqnarray*}
\end{theo}

\subsection{}
Letting $\cM=\cO_\cX$ in theorem \ref{thmloc}, which is justified by \cite[th\ptbl 5.4.3]{MSaito86}, one gets
\begin{cor}
If $f:\cX\rightarrow \CC$ is projective in a neighbourhood of $f^{-1}(0)$, one has for all $i\in\ZZ$
\begin{eqnarray*}
\dim\bH^{i-1}\lefpar f^{-1}(0),\phi_f\CC_\cX\rigpar &=&\dim\bH^i\lefpar f^{-1}(0),(\Omega_\cX^{\bbullet},df\wedge)\rigpar.
\end{eqnarray*}
\end{cor}

Notice that the corollary is well-known if one only assumes that, in a neighbourhood of $f^{-1}(0)$, the function $f$ has only a finite number of critical points; one then knows \cite{KSaito76} that the complex $(\Omega_{\cX}^{\bbullet},df\wedge)$ has cohomology in degree $\dim\cX$ at most and has support in the set of critical points; so is then the perverse sheaf $\phip_f\pCC_\cX$ where $\pCC_\cX\defin\CC_\cX[\dim\cX]$; in this case the corollary is consequence of the formula
\begin{eqnarray*}
\mu(f,x^o)&=&\dim\cO_{\cX,x^o}\left/\lefpar \frac{\partial f}{\partial x_1},\ldots ,\frac{\partial f}{\partial x_n}\rigpar\right.
\end{eqnarray*}
where $\mu(f,x^o)$ is the number of vanishing cycles of $f$ at $x^o$.

\subsection{}
Let \label{numlogan}$D$ be a divisor with normal crossings in $\cX$ and $j:\cU=\cX-D\hookrightarrow \cX$ the inclusion. Taking $\cM=\cO_\cX[*D]$ (\cf. \cite[\T2.d]{MSaito87}) one gets (see \T\ref{numrem2}):
\begin{cor}
Under the same assumptions as above one has, for all $i\in\ZZ$,
\begin{eqnarray*}
\dim\bH^{i-1}\lefpar f^{-1}(0),\phi_f\bR j_*\CC_\cU\rigpar &=&\dim\bH^i\lefpar f^{-1}(0),(\Omega_\cX^{\bbullet}(\log D),df\wedge) \rigpar.
\end{eqnarray*}
\end{cor}

\subsection{}
The proofs of theorems \ref{thmglob} and \ref{thmloc} are analogous: in the first case this is an exercise on filtered Fourier transform and in the second one an exercise on filtered microlocalization, once the results of \cite{MSaito86,MSaito87}, which are the difficult ones, are known. One proves in both cases the freeness of the analogues of the classical Brieskorn lattice, from which the theorems follow immediately.

In an appendix we give the proofs of theorems \ref{thmglob} and \ref{thmloc} suggested by the referee. These proofs are simpler than the proofs given below, as they do not use Fourier transform nor microlocalization.

\bigskip
\begin{small}
I thank F. Loeser for letting me know his notes of a talk by M. Kontsevitch and having raised my attention to this question. I thank the referee for his careful reading of the manuscript.
\end{small}

\section{Proof of theorem \protect\ref{thmglob}}
We keep notation of \T\ref{numnot}.

Let $(M,F)$ be a well-filtered coherent $\cD_X$-module. Let $\tau$ be a new variable. One puts on $M\otimes_\CC\CC[\tau,\tau^{-1}]$ the filtration
\begin{eqnarray*}
G_k\lefpar M\otimes_\CC\CC[\tau,\tau^{-1}]\rigpar &=&\ooplus_j F_{j+k}M \tau^{-j}
\end{eqnarray*}
so that $G_k=\tau^{k}G_0$ and
$$
\gr_0^G(M\otimes_\CC\CC[\tau,\tau^{-1}])=G_0/\tau^{-1}G_0\simeq\gr^F\!\!M.
$$
The proof of the theorem relies on the
\begin{prop}
Assume that the cohomology of the direct image $f_+M$ is holonomic and regular at infinity, and that $f_+(M,F)$ is strict. Then, for all $i\in\ZZ$, the $\CC[\tau^{-1}]$-module
$$
\bH^i\lefpar X,(\Omega_{X}^{\bullet}\ootimes_{\cO_X}G_0(M[\tau,\tau^{-1}]),\tau^{-1}\nabla-df\wedge)\rigpar
$$
is free of finite rank.
\end{prop}

\subsection{Computation of $\dim\bH^i\lefpar X,(\Omega_{X}^{\bbullet}\otimes_{\cO_X}M,\nabla-df\wedge)\rigpar$}\label{num21}
Consider the direct image $f_+M$. Then the cohomology modules $\cH^if_+M$ of $f_+M$ are regular holonomic (even at infinity) on the Weyl algebra $\CC[t]\langle\partial_t\rangle$, if $t$ denotes the coordinate on the affine line (see \eg. \cite{Borelal87}). The Fourier transform $\wwwh{\cH^if_+M}$ obtained by ``doing $\tau=\partial _t$ and $\partial_\tau=-t$'' is therefore a holonomic $\CC[\tau]\langle\partial_\tau\rangle$-module with a regular singularity at $\tau=0$ and possibly an irregular singularity at $\tau=\infty$. It has no other singularity, so that $\CC[\tau,\tau^{-1}]\otimes_{\CC[\tau]}\wwwh{\cH^if_+M}$ is a free $\CC[\tau,\tau^{-1}]$-module of rank $\mu_i$. Hence this number $\mu_i$ is equal to $\dim\wwwh{\cH^if_+M}/(\tau-1)\wwwh{\cH^if_+M}$ (one proves this first for $\CC[t]\langle\partial_t\rangle$-modules of the form $\CC[t]\langle\partial_t\rangle/(P)$ with $P\in\CC[t]\langle\partial_t\rangle$ nonzero and regular even at infinity, then, by an extension argument, for any regular holonomic $\CC[t]\langle\partial_t\rangle$-module; see for instance \cite[Chap\ptbl V]{Malgrange91}).

In order to express this dimension, one first computes the Fourier transform $\wwh{f_+M}$ of the complex $f_+M$. Let $\Afuc$ denote the affine line with coordinate $\tau$ and let $\pch:X\times \Afuc\rightarrow X$ and $\qc:X\times \Afuc\rightarrow \Afuc$ be the projections. Let $M[\tau]e^{-\tau f}$ be the $\cD_{X\times \Afuc}$-module $\pch^+M$ on which the action of $\cD$ is ``twisted'' by $e^{-\tau f}$. One then has
\begin{eqnarray*}
\wwh{f_+M}&=&\qc_+\lefpar M[\tau]e^{-\tau f}\rigpar
\end{eqnarray*}
(in order to understand better this formula, it may be useful to consider the direct image $i_+M=M[\partial_t]$ of $M$ by the inclusion $i:X\hookrightarrow X\times \AA^1$ defined by the graph of $f$ and to perform the partial Fourier transform on $i_+M$ with respect ot the $t$ variable).

Moreover, one knows that $\wwwh{\cH^if_+M}=\cH^i\wwh{f_+M}$ (see for instance \cite[App\ptbl 2]{Malgrange91}). One deduces from the freeness statement above that
\begin{eqnarray*}
\wwwh{\cH^if_+M}/(\tau-1)\wwwh{\cH^if_+M}&=&H^i\qc_+(Me^{-f}),
\end{eqnarray*}
$\qc:X\rightarrow \rm pt$ still denoting the restriction to $\tau=1$ of $\qc$. In this way one finds that
\begin{eqnarray*}
\mu_i&=&\dim\bH^i\lefpar X,(\Omega_{X}^{\bbullet+\dim X}\otimes_{\cO_X}M,\nabla-df\wedge)\rigpar.
\end{eqnarray*}

On the other hand, it is also known (see \cite[cor\ptbl 8.3]{Brylinski86} or \cite[Chap\ptbl V, prop\ptbl 1.5]{Malgrange91}) that
\begin{eqnarray*}
\mu_i&=&\sum_c\dim\phip_{t-c}\pDRa\lefpar \cH^if_+M\rigpar,
\end{eqnarray*}
where the sum is taken at singular points $c$ of $\cH^if_+M$ (at other points $c$, the corresponding term is zero), as asserted in remark \ref{remphi}.

As the functor $\phip$ preserves perversity (see \cite[cor\ptbl 1.7]{Brylinski86}) and commutes with the proper direct image, and as $\bR f_*$ commutes with $\pDRa$ in the present situation, one has
\begin{eqnarray*}
\mu_i&=&\sum_c\dim\cH^i(\phip_{t-c}\pDRa f_+\cM)\\
&=&\sum_c\dim\cH^i(\phip_{t-c}\bR f_*\pDRa \cM)\\
&=&\sum_c\dim \bH^i\lefpar f^{-1}(c),\phip_{f-c}\pDRa\cM\rigpar.
\end{eqnarray*}
One actually gets, as asserted in remark \ref{numrem},
\begin{eqnarray*}
\dim\bH^i\lefpar X,(\Omega_{X}^{\bbullet+\dim X} \otimes_{\cO_X}M,\nabla-df\wedge)\rigpar&=& \sum_{c\in\CC}\dim\bH^i\lefpar f^{-1}(c),\phip_{f-c}\pDRa(M)\rigpar,
\end{eqnarray*}
where all the terms in the sum except a finite number of them are zero.

\subsection{Proof of theorem \protect\ref{thmglob}}
The essential point is that mixed Hodge Modules satisfy the assumptions of the proposition in the situation of the theorem, as $f$ is projective (see \cite[th\ptbl 2.14]{MSaito87}), so that we will use it to prove the theorem.

We have seen that
\begin{eqnarray*}
\CC[\tau,\tau^{-1}]\ootimes_{\CC[\tau]}\wwwh{\cH^if_+M}&=& \bH^i\lefpar X,(\Omega_{X}^{\bullet+\dim X}\ootimes_{\cO_X}M[\tau,\tau^{-1}],\nabla-\tau df\wedge)\rigpar\\
&=&\bH^i\lefpar X,(\Omega_{X}^{\bullet+\dim X}\ootimes_{\cO_X}M[\tau,\tau^{-1}],\tau^{-1}\nabla-df\wedge)\rigpar
\end{eqnarray*}
is a free $\CC[\tau,\tau^{-1}]$-module of rank $\mu_i$. On the other hand, as the $\CC[\tau^{-1}]$-module
$$
\bH^i\lefpar X,(\Omega_{X}^{\bullet+\dim X}\ootimes_{\cO_X}G_0(M[\tau,\tau^{-1}]),\tau^{-1}\nabla-df\wedge)\rigpar
$$
is free, it is contained (by the natural map) in $\CC[\tau,\tau^{-1}]\ootimes_{\CC[\tau]}\wwwh{\cH^if_+M}$. Also, as it has finite rank and generates this module over $\CC[\tau,\tau^{-1}]$, it is a lattice in it, hence is $\CC[\tau^{-1}]$-module of rank $\mu_i$.

Since the terms of the complex $(\Omega_{X}^{\bullet+\dim X}\ootimes_{\cO_X}G_0(M[\tau,\tau^{-1}]),\tau^{-1}\nabla-df\wedge)$ as well as the hypercohomology modules are free $\CC[\tau^{-1}]$-modules, one gets, tensorizing by $\CC[\tau^{-1}]/\tau^{-1}\CC[\tau^{-1}]$,
\begin{eqnarray*}
\mu_i&=&\dim\bH^i\lefpar X,(\Omega_{X}^{\bullet+\dim X}\ootimes_{\cO_X}\gr^F\!\!M,\gr^F\!\nabla-df\wedge)\rigpar.\quad\qed
\end{eqnarray*}

\subsection{Filtered direct image}\label{numimdir}
Let $(M,F)$ be a well-filtered coherent $\cD_X$-module and let $f:X\rightarrow \AA^1$ be a regular function, which is proper on the support of $M$. If $\cD_{\AA^1\leftarrow X}$ denotes the transfer module associated to $f$, one puts $f_+M=\bR f_*\lefpar \cD_{\AA^1\leftarrow X}\ootimes_{\cD_X}^{\bL}M\rigpar$. Let $\omega_X=\Omega_{X}^{\dim X}$ be the sheaf of differential forms of maximum degree, equipped with its natural structure of right $\cD_X$-module and denote by $t$ the coordinate on $\AA^1$. Then $\cD_{\AA^1\leftarrow X}=\omega_X[\partial_t]$ as a $\cO_X$-module. The structure of right $\cD_X$-module on $\omega_X[\partial_t]$ is given by
\begin{eqnarray*}
\lefpar \sum_{j\geq 0}\omega_j\partial_{t}^{j}\rigpar \cdot \partial_{x_i}&=&\sum_{j\geq 0}\lefpar \omega_j\cdot\partial_{x_i}-\frac{\partial f}{\partial x_i}\omega_{j-1}\rigpar \partial _{t}^{j}
\end{eqnarray*}
and the structure of left $f^{-1}\cD_{\AA^1}$-module is given by
$$
\partial_t\cdot\lefpar \sum_{j\geq 0}\omega_j\partial_{t}^{j}\rigpar=\lefpar \sum_{j\geq 0}\omega_j\partial_{t}^{j+1}\rigpar\qqbox{and}t\cdot\lefpar \sum_{j\geq 0}\omega_j\partial_{t}^{j}\rigpar=-\lefpar \sum_{j\geq 1}j\omega_j\partial_{t}^{j-1}\rigpar.
$$

Let $\cR_FM\defin\oplus _{k\in\ZZ}F_k M\hbar^k\subset M[\hbar,\hbar^{-1}]$ be the Rees module on the new variable $\hbar$ associated to the filtration $F$. As $F$ is good, this is a coherent graded $\cR_F\cD_X$-module, $F\cD_X$ denoting the filtration by the order of differential operators. Considering the filtration $F\cD_{\AA^1\leftarrow X}$ defined by $F_k\cD_{\AA^1\leftarrow X}=\oplus _{j\leq k}\omega_X\partial_t^j$, one defines in the same way the filtered direct image by the formula
\begin{eqnarray*}
f_+\cR_FM&\defin&\bR f_*\lefpar \cR_F\cD_{\AA^1\leftarrow X}\ootimes_{\cR_F\cD_X}^{\bL}\cR_FM\rigpar,
\end{eqnarray*}
which is an object of $D_{\rm coh}^{b}(\cR_F\cD_{\AA^1})$ (see for instance \cite[\S\T2.1, 2.3]{MSaito86}).

\bigskip
The complex
$$\lefpar \Omega_{X}^{\bbullet+\dim X}\otimes_{\cO_X}\cD_X[\partial _t],\nabla-\partial _t df\wedge\rigpar
$$
with $\nabla(\omega\otimes P)=d\omega\otimes P+\sum_i dx_i\wedge\omega\otimes \partial_{x_i}P$, is a free resolution of the $(f^{-1}\cD_{\AA^1},\cD_X)$-bimodule $\cD_{\AA^1\leftarrow X}$. One verifies that it is strict for the diagonal filtration $G_k\cD_X[\partial _t]=\sum_{i+j=k}F_i\cD_X\partial _t^j$, that is, the complex
$$\lefpar \Omega_{X}^{\bbullet+\dim X}\otimes_{\cO_X}\cR_{G[-\bbullet]}\cD_X[\partial _t],\nabla-\partial _t df\wedge\rigpar
$$
is a resolution of $\cR_F\cD_{\AA^1\leftarrow X}$.

\bigskip
For $(M,F)$ as above, denote by $G_k(M[\partial _t])=\sum_{i+j=k}F_i M\partial_{t}^{j}$. One then has
\begin{eqnarray*}
f_+\cR_FM&=&\bR f_*\lefpar \Omega_{X}^{\bbullet+\dim X}\otimes_{\cO_X}\cR_{G[-\bbullet]}M[\partial _t],\nabla-\partial _t df\wedge\rigpar.
\end{eqnarray*}

\begin{lem}
For $i\in\ZZ$, the following conditions are equivalent:
\begin{enumerate}
\item
The cohomology modules $\cH^if_+\cR_FM$ are $\cR_F\cD_{\AA^1}$-modules without $\hbar$-torsion.
\item
For all $k\in\ZZ$, the following natural morphism is injective:
\begin{multline*}
\cH^i\bR f_*\lefpar \Omega_{X}^{\bbullet+\dim X}\otimes_{\cO_X}G_{k+\bbullet}M[\partial _t],\nabla-\partial _t df\wedge\rigpar\\
\longrightarrow \cH^i\bR f_*\lefpar \Omega_{X}^{\bbullet+\dim X}\otimes_{\cO_X}M[\partial _t],\nabla-\partial _t df\wedge\rigpar=\cH^if_+M.
\end{multline*}
\end{enumerate}
\end{lem}

\begin{proof}
As $f$ is proper on $\supp M$, $\cH^i\bR f_*$ commutes with direct limit and the second condition is equivalent to the injectivity of the natural map associated to localization with respect to $\hbar$:
\begin{eqnarray*}
\cH^i f_+\cR_FM&\longrightarrow &\cH^i f_+\lefpar \cR_FM\ootimes_{\CC[\hbar]}\CC[\hbar,\hbar^{-1}]\rigpar.
\end{eqnarray*}
The right-hand term is precisely the localized module $\lefpar\cH^i f_+ \cR_FM\rigpar \ootimes_{\CC[\hbar]}\CC[\hbar,\hbar^{-1}]$, as $\CC[\hbar,\hbar^{-1}]$ is flat on $\CC[\hbar]$, and the injectivity of this map is then equivalent to the $\hbar$-torsionlessness of $\cH^i f_+ \cR_FM$.
\end{proof}

Generally speaking, the modules $\cH^i f_+ M$ come equipped with a natural filtration denoted by $F\cH^i f_+ M$: it is defined as the image of the map introduced in the statement 2 of the lemma. One then has
\begin{eqnarray*}
\cR_F\cH^i f_+ M&=&\cH^i f_+ \cR_FM/\hbar\mbox{-torsion}.
\end{eqnarray*}
So the filtration $F\cH^i f_+ M$ is good.

When the conditions of the lemma are satisfied for all $i\in\ZZ$, one says that $f_+(M,F)$ is {\em strict}. It is equivalent to saying that one has for all $i$
\begin{eqnarray*}
\cR_F\cH^i f_+ M&=&\cH^i f_+ \cR_FM.
\end{eqnarray*}

\subsection{Proof of the proposition}
We will identify $\cD_{\AA^1}$-modules and $\CC[t]\langle\partial_t\rangle$-modules {\em via} the global sections functor $\Gamma(\AA^1,\bullet)$. Let $F\CC[t]\langle\partial_t\rangle$ and $F\CC[t]\langle\partial_t,\partial_{t}^{-1}\rangle$ be the filtrations by the degree in $\partial_t$. Then $\cR_F\CC[t]\langle\partial_t,\partial_{t}^{-1}\rangle$ is flat on $\cR_F\CC[t]\langle\partial_t\rangle$ as a graded module (the proof is analogous to that of the lemma in \T\ref{nummicro} below).

If the cohomology of $f_+\cR_FM$ has no $\hbar$-torsion (\ie. if $f_+(M,F)$ is strict), then so is that of
$$
\cR_F\CC[t]\langle\partial_t,\partial_{t}^{-1}\rangle\ootimes_{\cR_F\CC[t]\langle\partial_t\rangle}\Gamma(\AA^1,f_+\cR_FM)
$$
because of flatness. Let us compute this complex: because of flatness and using the projection formula, one has
\begin{multline*}
\cR_F\CC[t]\langle\partial_t,\partial_{t}^{-1}\rangle\ootimes_{\cR_F\CC[t]\langle\partial_t\rangle}\Gamma(\AA^1,f_+\cR_FM)\\
\shoveleft =\bR\Gamma \lefpar \AA^1,\bR f_*\lefpar f^{-1}\cR_F\cD_{\AA^1}\langle\partial_{t}^{-1}\rangle\ootimes_{\cR_F\cD_{\AA^1}}(\cR_F\cD_{\AA^1\leftarrow X}\ootimes_{f^{-1}\cR_F\cD_X}^{\bL}\cR_FM)\rigpar\rigpar\\
=\bR\Gamma\lefpar X,(\Omega_{X}^{\bbullet+\dim X}\ootimes_{\cO_X}\cR_{G[-\bbullet]}(M[\partial_t,\partial_{t}^{-1}]),\nabla-\partial_t df\wedge)\rigpar.
\end{multline*}
As in lemma of \T\ref{numimdir}, owing to $G_k(M[\partial_t,\partial_{t}^{-1}])=\partial_{t}^{k}G_0(M[\partial_t,\partial_{t}^{-1}])$, the absence of $\hbar$-torsion is equivalent to the injectivity, for all $i\in\ZZ$, of the map
\begin{multline*}
\bH^i\lefpar X,(\Omega_{X}^{\bbullet+\dim X}\ootimes_{\cO_X}G_{-\bbullet}(M[\partial_t,\partial_{t}^{-1}]),\nabla-\partial_t df\wedge)\rigpar\\
\longrightarrow \bH^i\lefpar X,(\Omega_{X}^{\bbullet+\dim X}\ootimes_{\cO_X}M[\partial_t,\partial_{t}^{-1}],\nabla-\partial_t df\wedge)\rigpar.
\end{multline*}
The right-hand side, which can be identified to $\CC[\partial_t,\partial_{t}^{-1}]\otimes_{\CC[\partial_t]}\wwwh{\cH^if_+M}$, is $\CC[\partial_t,\partial_{t}^{-1}]$-free of rank $\mu_i$ if $\cH^if_+M$ is holonomic and regular at infinity.

On the other hand, the image of the map above can be obtained from the good filtration $F\cH^if_+M$: it is equal to $\sum_{j\in\ZZ}\partial_{t}^{-j}\alpha\lefpar F_{i+j-\dim X} \cH^if_+M\rigpar$, where $\alpha:\cH^if_+M\rightarrow \cH^if_+M[\partial_{t}^{-1}]$ denotes the localization map. This image is then a free $\CC[\partial_{t}^{-1}]$-module of rank $\mu_i$ which generates $\CC[\partial_t,\partial_{t}^{-1}]\otimes_{\CC[\partial_t]}\wwwh{\cH^if_+M}$ over $\CC[\partial_t,\partial_{t}^{-1}]$ (see for instance \cite[prop\ptbl 2.1]{Bibi96b}): by an extension argument, it is enough to prove this for any well-filtered holonomic $\CC[t]\langle\partial_t\rangle$-module $(N,N_\bbullet)$ with $N_j=0$ for $j<0$ and $N_j=\partial_t^jN_0$ for $j\geq 0$, such that $N$ is regular at infinity; for the same reason, it is enough to consider $N=\CC[t]\langle\partial_t\rangle/(P)$ with $P=\sum_{i=0}^{d}\partial_{t}^{i}a_i(t)$, $a_i\in\CC[t]$, $a_d\not\equiv0$ and $\deg a_i<\deg a_d$ for $i<d$, with the filtration induced by the degree in $\partial_t$; in this case the result is easy.

As the map is injective, its source verifies the same properties. Finally, one has
\begin{multline*}
\bH^i\lefpar X,(\Omega_{X}^{\bbullet+\dim X}\ootimes_{\cO_X}G_{-\bbullet}(M[\partial_t,\partial_{t}^{-1}]),\nabla-\partial_t df\wedge)\rigpar\\
=\bH^i\lefpar X,(\Omega_{X}^{\bbullet+\dim X}\ootimes_{\cO_X}G_0(M[\partial_t,\partial_{t}^{-1}]),\partial_{t}^{-1}\nabla-df\wedge)\rigpar.\quad\qed
\end{multline*}

\section{Proof of theorem \protect\ref{thmloc}}
The results of \T\ref{numimdir} can be straightforwardly adapted to the holomorphic case.

\subsection{Microlocalization}\label{nummicro}
Let $\wh\cE_\CC$ be the restriction to $\CC=\CC\times {\bf1}\subset T^*\CC$ of the sheaf of formal microdifferential operators: a section of $\wh\cE_\CC$ on $U\subset\CC$ is a formal series $\sum_{n\leq n_0}a_n(t)\partial_{t}^{n}$ where the $a_n$ are holomorphic functions on $U$. This sheaf is filtered by
$$
F_k\wh\cE_\CC\defin\wh\cE_\CC(k)=\wh\cE_\CC(0)\partial_{t}^{k}= \partial_{t}^{k}\wh\cE_\CC(0),
$$
where $\wh\cE_\CC(0)$ is made of formal series $\sum_{n\leq 0}a_n(t)\partial_{t}^{n}$. One has $\wh\cE_\CC(0)/\wh\cE_\CC(-1)=\cO_\CC$ and $\gr^F\!\wh\cE_\CC=\cO_\CC[\wt\tau,\wt\tau^{-1}]$ if $\wt\tau$ denotes the class of $\partial_t$ in $\wh\cE_\CC(1)/\wh\cE_\CC(0)$. Let $\cR_F\wh\cE_\CC=\oplus_k\wh\cE_\CC(k)\hbar^k$ be the Rees ring associated to this filtration.

\begin{lem}
The ring $\cR_F\wh\cE_\CC$ is flat on $\cR_F\cD_\CC$ as a graded module.
\end{lem}

\begin{proof}
The question is local. Denote $\cD$ and $\wh\cE$ the germs at $c\in\CC$ of $\cD_\CC$ and $\wh\cE_\CC$. Let $L$ and $L'$ be free $\cD$-modules equipped with filtrations of the kind $(L,F)=\oplus _k(\cD,F[n_k])$, $(L',F)=\oplus _\ell(\cD,F[n'_\ell])$. Let $\varphi:(L',F)\rightarrow (L,F)$ be a strict morphism, \ie. such that $\im\varphi\cap F_\bbullet L=\varphi(F_\bbullet L')$. It is a matter of verifying that $1\otimes \varphi:\wh\cE\otimes_\cD(L',F)\rightarrow \wh\cE\otimes_\cD(L,F)$ is still strict: indeed, if this assertion is proved, consider an ideal $\cR_F I$ of $\cR_F\cD$, where $(I,F)$ is a well-filtered coherent ideal of $\cD$; there exists a presentation of the kind
$$
\cR_FL'\MRE{\cR_F\varphi}\cR_FL\MRE{}\cR_F I\subset\cR_F\cD
$$
and one wants to show that the sequence obtained after $\cR_F\wh\cE\otimes_{\cR_F\cD}\bullet$ remains exact; it is thus a matter of seeing that the sequence $(\wh\cE\otimes_\cD L',F)\rightarrow (\wh\cE\otimes_\cD L,F)\rightarrow (\wh\cE,F)$ is exact and strict; exactness follows from flatness of $\wh\cE$ over $\cD$ (see for instance \cite{Bjork79,Schapira85}) and strictness is then the assertion mentioned above, that we now show.

In the bases of $L$ and $L'$ that we have considered, the morphism $\varphi$ is expressed as right multiplication by a matrix $A(t,\partial _t)$ with elements in $\cD$. Let $a'=\sum_{i\leq i_0}a'_i(t)\partial_{t}^{i}\in F_{p'}(\wh\cE\otimes L')=\oplus_\ell \wh\cE(p'-n'_\ell)$ and assume that $1\otimes \varphi(a')\in F_p(\wh\cE\otimes L)=\oplus_k\wh\cE(p-n_k)$. We want to find $a''\in F_p(\wh\cE\otimes L')$ such that $1\otimes \varphi(a'')=1\otimes \varphi(a')$.

We may choose $i_1\leq i_0$ so that if one sets $a'=b'+c'$ with $b'=\sum_{i\leq i_1}a'_i(t)\partial_{t}^{i}$, one has $b'\in F_p(\wh\cE\otimes L')$ and $1\otimes \varphi(b')\in F_p(\wh\cE\otimes L)$.

Let $m\in\NN$ be such that $\partial_{t}^{m}c'\in L'$. Then $\varphi(\partial_{t}^{m}c')\in L\cap F_{p+m}(\wh\cE\otimes L)=F_{p+m}L$. As $\varphi$ is strict, there exists $c''\in F_{p+m}L'$ such that $\varphi(c'')=\varphi(\partial_{t}^{m}c')$. We may then take $a''\defin b'+\partial_{t}^{-m}c''$.
\end{proof}

Let $(\cM,F)$ and $f:\cX\rightarrow \CC$ be as in \T\ref{numimdir}. One may microlocalize the filtered direct image $f_+\cR_F\cM$ by putting
\begin{eqnarray*}
(f_+\cR_F\cM)^{\!^\mu}&=&\cR_F\wh\cE_\CC\ootimes_{\cR_F\cD_\CC}f_+\cR_F\cM.
\end{eqnarray*}
The projection formula and the lemma above show that
\begin{eqnarray*}
(f_+\cR_F\cM)^{\!^\mu}&=&\bR f_*\lefpar f^{-1}\cR_F\wh\cE_\CC\ootimes_{f^{-1}\cR_F\cD_\CC}\lefpar \cR_F\cD_{\CC\leftarrow\cX}\ootimes_{\cR_F\cD_\CC}^{\bL}\cR_F\cM\rigpar \rigpar.
\end{eqnarray*}

As $f^{-1}\cR_F\wh\cE_\CC\otimes_{f^{-1}\cR_F\cD_\CC}\cR_G(\cM[\partial_t])$ has no $\hbar$-torsion (see the flatness lemma above), there exists a unique filtration $G_\bbullet\lefpar f^{-1}\wh\cE_\CC\otimes_{f^{-1}\cD_\CC}\cM[\partial_t]\rigpar$ such that
\begin{eqnarray*}
\cR_G\lefpar f^{-1}\wh\cE_\CC\otimes_{f^{-1}\cD_\CC}\cM[\partial_t]\rigpar &=&f^{-1}\cR_F\wh\cE_\CC\otimes_{f^{-1}\cR_F\cD_\CC}\cR_G(\cM[\partial_t]).
\end{eqnarray*}
This filtration is defined by the formula
\begin{eqnarray*}
G_k\lefpar f^{-1}\wh\cE_\CC\otimes_{f^{-1}\cD_\CC}\cM[\partial_t]\rigpar&=&\sum_{i+j=k}\image\lefpar f^{-1}\wh\cE_\CC(i)\otimes_{f^{-1}\cO_\CC}G_j(\cM[\partial_t])\rigpar,
\end{eqnarray*}
which shows that for all $k\in\ZZ$ one has
\begin{eqnarray*}
G_k\lefpar f^{-1}\wh\cE_\CC\otimes_{f^{-1}\cD_\CC}\cM[\partial_t]\rigpar&=&\partial_{t}^{k}G_0\lefpar f^{-1}\wh\cE_\CC\otimes_{f^{-1}\cD_\CC}\cM[\partial_t]\rigpar.
\end{eqnarray*}

On the other hand, the flatness lemma above shows that the exact sequence
$$
0\longrightarrow \cR_G(\cM[\partial_t])\MRE{\hbar} \cR_G(\cM[\partial_t])\MRE{} \gr^G\!(\cM[\partial_t])\longrightarrow 0
$$
remains exact after tensorizing by $f^{-1}\cR_F\wh\cE_\CC$, so that one has
\begin{eqnarray*}
\gr^G\!\lefpar f^{-1}\wh\cE_\CC\ootimes_{f^{-1}\cD_\CC}\cM[\partial_t]\rigpar&=&f^{-1}\gr^F\!\wh\cE_\CC\ootimes_{f^{-1}\gr^F\!\cD_\CC}\gr^G\!(\cM[\partial_t])\\
&=&f^{-1}\cO_\CC[\wt\tau,\wt\tau^{-1}]\ootimes_{f^{-1}\cO_\CC[\wt\tau]}(\gr^F\!\!\cM)[\wt\tau]\\
&=&\gr^F\!\!\cM[\wt\tau,\wt\tau^{-1}],
\end{eqnarray*}
where the grading of the right-hand side is the diagonal grading. One deduces in particular that for all $k\in\ZZ$ one has
$$
\gr_k^G\lefpar f^{-1}\wh\cE_\CC\ootimes_{f^{-1}\cD_\CC}\cM[\partial_t]\rigpar\arrowl{\wt\tau^k}{\sim}\gr_0^G\lefpar f^{-1}\wh\cE_\CC\ootimes_{f^{-1}\cD_\CC}\cM[\partial_t]\rigpar \MRE{\sim} \gr^F\!\!\cM.
$$
One finally gets
\begin{eqnarray*}
(f_+\cR_F\cM)^{\!^\mu}&=&\bR f_*\lefpar \Omega_{\cX}^{\bbullet+\dim\cX}\ootimes_{\cO_\cX}\cR_{G[-\bbullet]}(f^{-1}\wh\cE_\CC\ootimes_{f^{-1}\cD_\CC}\cM[\partial_t]), \nabla-\partial_t df\wedge\rigpar
\end{eqnarray*}
and
\begin{multline*}
\bR f_*\gr_0^G\lefpar \Omega_{\cX}^{\bbullet+\dim\cX}\ootimes_{\cO_\cX}(f^{-1}\wh\cE_\CC\ootimes_{f^{-1}\cD_\CC}\cM[\partial_t]), \nabla-\partial_t df\wedge\rigpar\\
= \bR f_*\lefpar \Omega_{\cX}^{\bbullet+\dim\cX}\ootimes_{\cO_\cX}\gr^F\!\!\cM,\gr^F\!\nabla-df\wedge\rigpar.
\end{multline*}

\subsection{Microlocal degeneracy at $E_1$}
In order to simplify the notation, we will set in the following $\theta=\partial_{t}^{-1}$.

\begin{lem}
For all $i\in\ZZ$ the following conditions are equivalent:
\begin{enumerate}
\item
$\cH^i\bR f_*G_0\lefpar \Omega_{\cX}^{\bbullet+\dim\cX}\ootimes_{\cO_\cX}(f^{-1}\wh\cE_\CC\ootimes_{f^{-1}\cD_\CC}\cM[\partial_t]), \nabla-\partial_t df\wedge\rigpar$ has no $\theta$-torsion.
\item
The following natural morphism is injective:
\begin{eqnarray*}
\cH^i\bR f_*G_0\lefpar \Omega_{\cX}^{\bbullet+\dim\cX}\ootimes_{\cO_\cX}(f^{-1}\wh\cE_\CC\ootimes_{f^{-1}\cD_\CC}\cM[\partial_t]), \nabla-\partial_t df\wedge\rigpar&\longrightarrow &\cH^i(f_+\cM)^{\!^\mu}.
\end{eqnarray*}
\item
$\cH^i(f_+\cR_F\cM)^{\!^\mu}$ has no $\hbar$-torsion.
\end{enumerate}
Moreover, these conditions are satisfied for all $i$ if $f_+(\cM,F)$ is strict.
\end{lem}

\begin{proof}
One has for all $k\in\ZZ$
\begin{eqnarray*}
f^{-1}\wh\cE_\CC\ootimes_{f^{-1}\cD_\CC}\cM[\partial_t]&=&f^{-1}\wh\cE_\CC\ootimes_{f^{-1}\wh\cE_\CC(0)}G_0\lefpar f^{-1}\wh\cE_\CC\ootimes_{f^{-1}\cD_\CC}\cM[\partial_t]\rigpar\\
&=&\CC\lcr\theta\rcr[\theta^{-1}]\ootimes_{\CC\lcr\theta\rcr}G_0\lefpar f^{-1}\wh\cE_\CC\ootimes_{f^{-1}\cD_\CC}\cM[\partial_t]\rigpar.
\end{eqnarray*}
Denote by $G_0$ the complex which appears at point 1. The map of point 2 may thus be identified with
\begin{eqnarray*}
\cH^i\bR f_*G_0&\longrightarrow &\CC\lcr\theta\rcr[\theta^{-1}]\ootimes_{\CC\lcr\theta\rcr}\cH^i\bR f_*G_0.
\end{eqnarray*}
The equivalence of the statements 1 and 2 is then clear.

The equivalence of the statements 2 and 3 can be shown as in the lemma of \T\ref{numimdir}. Finally, if $f_+(\cM,F)$ is strict, $\cH^if_+\cR_F\cM$ has no $\hbar$-torsion by definition (for all $i\in\ZZ$), and from the flatness lemma one deduces that neither has $\cH^i(f_+\cR_F\cM)^{\!^\mu}$, hence the statement 3 is verified.
\end{proof}

\subsection{The holonomic case}\label{numcashol}
If $\cH^if_+\cM$ is holonomic (for instance if $\cM$ is holonomic), the microlocal module $\cH^i(f_+\cM)^{\!^\mu}$ has punctual support in $\CC$. We will assume in the following that this module has support $\{t=0\}$. As $f_+(\cR_F\cM)^{\!^\mu}$ has coherent cohomology over $\cR_F\wh\cE_\CC$, the image of the morphism of the statement 2 of the lemma is $\wh\cE_\CC(0)$-coherent, hence is free over $\CC\lcr\theta\rcr$, and of rank $\mu_i=\dim_{\CC\lcr\theta\rcr[\theta^{-1}]}\cH^i(f_+\cM)^{\!^\mu}$, according to the preparation theorem (see for instance \cite[\T3]{Malgrange81}). It follows then that, for a holonomic $\cM$, the conditions of the lemma are equivalent to
\begin{enumerate}
\item[1 bis.]
$\cH^i\bR f_*G_0\lefpar \Omega_{\cX}^{\bbullet+\dim\cX}\ootimes_{\cO_\cX}(f^{-1}\wh\cE_\CC\ootimes_{f^{-1}\cD_\CC}\cM[\partial_t]), \nabla-\partial_t df\wedge\rigpar$ is free over $\CC\lcr\theta\rcr$ of rank $\mu_i$.
\end{enumerate}

One concludes

\begin{lem}
If $f_+\cM$ has holonomic cohomology and $f_+(\cM,F)$ is strict, one has
\begin{eqnarray*}
\dim \bH^i\lefpar f^{-1}(0),\lefpar \Omega_{\cX}^{\bbullet+\dim\cX}\ootimes_{\cO_\cX}\gr^F\!\!\cM,\gr^F\!\nabla-df\wedge\rigpar\rigpar &=&\dim \bH^i\lefpar f^{-1}(0),\phip_f\pDRa\cM\rigpar.
\end{eqnarray*}
\end{lem}

\begin{proof}
Indeed, as the cohomology of the complex $\bR f_* G_0$ is $\CC\lcr\theta\rcr$-free and that the terms of the complex $G_0$ have no $\theta$-torsion, one has
\begin{eqnarray*}
\cH^i\bR f_*(G_0/\theta G_0)&=&G_0\cH^i(f_+\cM)^{\!^\mu}/\theta G_0\cH^i(f_+\cM)^{\!^\mu}.
\end{eqnarray*}
As we have seen above, the left-hand side can be identified to the left-hand side of the lemma, and the right-hand side has dimension $\mu_i$. It remains to identify $\mu_i$ with the right-hand side of the lemma. Recall that, according to the index theorem of Kashiwara (see \cite{Kashiwara70}, see also \cite[Chap\ptbl IV, cor\ptbl 4.2]{Malgrange91}), one has $\mu_i=\dim\phip_t\pDRa(\cH^if_+\cM)$. One may now conclude as in \T\ref{num21}.
\end{proof}

\begin{proof}[End of the proof]
Theorem \ref{thmloc} is now a consequence of the fact that mixed Hodge Modules satisfy the properties of the lemma above if $f$ is projective on $\supp\cM$ (see \cite[th\ptbl 2.14]{MSaito87}).
\end{proof}

\section{Some remarks}

\subsection{Interpretation of the graded complex}\label{numrem1}
Let $f:\cX\rightarrow \CC$ be as in \T\ref{numan}. Let $\sigma:\cX\hookrightarrow T^*\cX$ be the section defined by the $1$-form $df$. Let $\cO_\cX[T\cX]=\gr^F\!\cD_\cX$ be the sheaf of holomorphic functions on $T^*\cX$ which are polynomials in the fibres of $\pi:T^*\cX\rightarrow \cX$. Let $\omega_\cX=\Omega_{\cX}^{\dim\cX}$. Consider the $\cO_\cX[T\cX]$-module $\sigma_*\omega_\cX$: it coincides with $\omega_\cX$ as a $\cO_\cX$-module, when one considers $\cO_\cX$ as a subring of $\cO_\cX[T\cX]$ {\em via} $\pi^*$.

The Spencer resolution $(\Omega_{\cX}^{\bbullet+\dim\cX}\otimes_{\cO_\cX}\cD_\cX,\nabla)$ of $\omega_\cX$ as a right $\cD_\cX$-module, with, in local coordinates , $\nabla(\omega\otimes P)=\sum_i(dx_i\wedge\omega)\otimes \partial_{x_i}P$, is strictly filtered, when one equips $\Omega_{\cX}^{\ell+\dim\cX}\otimes_{\cO_\cX}\cD_\cX$ with the shifted filtration
\begin{eqnarray*}
F_k\lefpar \Omega_{\cX}^{\ell+\dim\cX}\otimes_{\cO_\cX}\cD_\cX\rigpar &\defin&\Omega_{\cX}^{\ell+\dim\cX}\otimes_{\cO_\cX}F_{k+\ell}\cD_\cX.
\end{eqnarray*}
After gradation one deduces a resolution of $\omega_\cX$ by locally free $\gr^F\!\cD_\cX$-modules:
\begin{eqnarray*}
(\Omega_{\cX}^{\bbullet+\dim\cX}\otimes_{\cO_\cX}\gr^F\!\cD_\cX,\gr^F\!\nabla)&\isom&\omega_\cX.
\end{eqnarray*}

In an analogous way one has
\begin{lem}
The complex $\lefpar \Omega_{\cX}^{\bbullet+\dim\cX}\otimes_{\cO_\cX}\gr^F\!\cD_\cX,\gr^F\!\nabla-df\wedge\rigpar $ is a $\gr^F\!\cD_\cX$-locally free resolution of $\sigma_*\omega_\cX$.
\end{lem}

\begin{proof}
In local coordinates $(x_1,\ldots ,x_n)$ on $\cX$, this complex can be identified with the Koszul complex associated to the regular sequence $\lefpar \dpl\xi_1-\frac{\partial f}{\partial x_1},\ldots ,\xi_n-\frac{\partial f}{\partial x_n}\rigpar$ of $\cO_\cX[\xi_1,\ldots ,\xi_n]$, hence the exactness.

On the other hand, $\sigma_*\omega_\cX$ is equal to the $\cO_\cX$-module $\omega_\cX$ on which $\xi_i$ acts as the product by $\partial f/\partial x_i$; this allows one to conclude.
\end{proof}

Let $(\cM,F)$ be a well-filtered coherent $\cD_\cX$-module. One deduces from the previous lemma that
\begin{eqnarray*}
\sigma_*\omega_\cX\ootimes_{\gr^F\!\cD_\cX}^{\bL}\gr^F\!\!\cM&=& \lefpar \Omega_{\cX}^{\bbullet+\dim\cX}\otimes_{\cO_\cX}\gr^F\!\!\cM,\gr^F\!\nabla-df\wedge\rigpar
\end{eqnarray*}
in $D_{\rm coh}^{b}(\cO_\cX)$. The cohomology sheaves of this complex have thus support in the inverse image by $\sigma$ of the characteristic variety $\car\cM=\supp\gr^F\!\!\cM$, in other words, in the intersection of $\car\cM$ with the image of the section $df$.

\bigskip
If $\cM$ is holonomic, the variety $\car\cM$ is Lagrangian. It is then well-known that, locally on $\cX$, $\sigma^{-1}(\car\cM)$ is contained in at most one fibre $f=\rm constant$.

In particular, if $f:\cX\rightarrow \CC$ is proper, there exists a neighbourhood of $f^{-1}(0)$ such that the complex $\lefpar \Omega_{\cX}^{\bbullet+\dim\cX}\otimes_{\cO_\cX}\gr^F\!\!\cM,\gr^F\!\nabla-df\wedge\rigpar$ has $\cO_\cX$-coherent cohomology with support contained in $f^{-1}(0)$.

If moreover $\cM$ is regular holonomic, one knows (see for instance \cite[p\ptbl 16]{Brylinski86}) that in a neighbourhood of any point $x^o\in\sigma^{-1}(\car\cM)$, the set $\sigma^{-1}(\car\cM)$ contains the support of the vanishing cycle complex $\phip_{f-f(x^o)}\pDRa\cM$.

\subsection{Logarithmic complexes}\label{numrem2}
Let $D$ be a divisor with normal crossings of a complex analytic manifold $\cX$. Equip the sheaf $\cO_\cX[*D]$ of meromorphic functions on $\cX$ with poles along $D$ with the increasing filtration $F_\bbullet$ by the pole order (\cf. \cite[p\ptbl 80]{Deligne70}). Put $\cO_\cX[*D][\tau]=\cO_\cX[*D]\otimes_\CC\CC[\tau]$ and set
\begin{eqnarray*}
G_k\cO_\cX[*D][\tau]&=&\sum_{i+j=k}F_i\cO_\cX[*D]\tau^j.
\end{eqnarray*}
One has $\gr^G\!\cO_\cX[*D][\tau]=\gr^F\!\cO_\cX[*D]\otimes_\CC\CC[\tau]$ equipped with the diagonal gradation. One gets a filtration $G_\bbullet(\Omega_{\cX}^{\bbullet}[*D][\tau],d-\tau df\wedge)$ as in \T\ref{numfilt}, using the fact that the differential $d-\tau df\wedge$ has degree $1$ with respect to $G_\bbullet$. The logarithmic sub-complex $(\Omega_{\cX}^{\bbullet}(\log D)[\tau],d-\tau df\wedge)$ inherits thus a filtration $G_\bbullet$.
\begin{lem}
The inclusion of filtered complexes
\begin{eqnarray*}
G_\bbullet(\Omega_{\cX}^{\bbullet}(\log D)[\tau],d-\tau df\wedge)&\hookrightarrow &G_\bbullet(\Omega_{\cX}^{\bbullet}[*D][\tau],d-\tau df\wedge)
\end{eqnarray*}
is a filtered quasi-isomorphism.
\end{lem}

\begin{proof}
It is a matter of verifying that the induced morphism on the graded complexes is a quasi-isomorphism. One has
\begin{eqnarray*}
\gr^G\!(\Omega_{\cX}^{\bbullet}(\log D)[\tau],d-\tau df\wedge)&=&(\Omega_{\cX}^{\bbullet}(\log D)[\tau],-\tau df\wedge).
\end{eqnarray*}
Consider then the filtration induced by $F_\bbullet$ on the complexes $\gr^G\!$. One shows that the $G$-graded morphism induces a $F$-filtered quasi-isomorphism: take thus the $F$-graded object and consider the graded morphism.

One has
\begin{eqnarray*}
\gr^F\!\gr^G\!(\Omega_{\cX}^{\bbullet}(\log D)[\tau],d-\tau df\wedge)&=&(\Omega_{\cX}^{\bbullet}(\log D)[\tau],0)
\end{eqnarray*}
and, identifying $\gr^F\!\gr^G\!\cO_\cX[*D][\tau]$ with $\gr^F\!\cO_\cX[*D][\tau]$, one has
\begin{eqnarray*}
\gr^F\!\gr^G\!(\Omega_{\cX}^{\bbullet}[*D][\tau],d-\tau df\wedge)&=&(\Omega_{\cX}^{\bbullet}\otimes_{\cO_\cX}\gr^F\!\cO_\cX[*D][\tau],\gr^F\!d)
\end{eqnarray*}
where $\gr^F\!d$ is seen as a graded operator of degree $1$ with respect to $F$. The quasi-isomorphism between both complexes now follows from \cite[II 3.13]{Deligne70}.
\end{proof}

\begin{proof}[Proof of corollaries \protect\ref{numlog} and \protect\ref{numlogan}]
Forget the filtration $G_\bbullet$ in the previous lemma. As the terms of the complexes are free over $\CC[\tau]$, one deduces from this lemma, after tensorizing with $\CC[\tau]/(\tau-1)$, a quasi-isomorphism
\begin{eqnarray*}
(\Omega_{\cX}^{\bbullet}(\log D),d-df\wedge)&\isom &(\Omega_{\cX}^{\bbullet}[*D],d-df\wedge).
\end{eqnarray*}

Consider now the quasi-isomorphism at the level of graded complexes given by the lemma. Then, ``putting $\tau=1$'' as above, one gets in the same way a quasi-isomorphism
\begin{eqnarray*}
(\Omega_{\cX}^{\bbullet}(\log D),-df\wedge)&\isom &(\Omega_{\cX}^{\bbullet}\ootimes_{\cO_\cX}\gr^F\!\cO_\cX[*D],\gr^F\!d-df\wedge)
\end{eqnarray*}
and the corollary \ref{numlogan} can be straightforwardly deduced from theorem \ref{thmloc}.

\medskip
The proof of corollary \ref{numlog} is identical, taking $\cX=\wt X^{\rm an}$ and exchanging above $\cO_\cX$ with $\cO_{\wt X^{\rm an}}[*D_\infty]$ with $D_\infty=\wt X-X$.
\end{proof}

\subsection{The K\"ahler case}
It follows from \cite{MSaito90} that theorems \ref{thmglob} and \ref{thmloc} remain true if one only assumes that $f$ is proper and K\"ahler on $\supp M$ or $\supp \cM$, imposing nevertheless that $M$ or $\cM$ has quasi-unipotent monodromy along any germ of holomorphic function $g:\Xt^{\rm an}\rightarrow \CC$ or $g:\cX\rightarrow \CC$, which is the case for mixed Hodge Modules with geometric origin, according to the monodromy theorem.

\subsection{Some questions}\mbox{ }\par
(1) Let $\crit(f,\cM)=\sigma^{-1}(\car\cM)$ be the critical set of $f$ with respect to $\cM$, where $\sigma$ denotes as above the section of $T^*\cX$ defined by $df$. Instead of assuming that $f$ is projective on $\supp\cM$, assume only that $f^{-1}(0)\cap\crit(f,\cM)$ is projective. Does theorem \ref{thmloc} remain valid under this assumption?

(2) Is it possible to give a purely algebraic proof, as in \cite{D-I87}, of corollaries \ref{numkont} and \ref{numlog}?

\section*{Appendix}
We give here direct proofs of theorems \ref{thmglob} and \ref{thmloc} suggested by the referee.

\begin{proof}[Proof of theorem \protect\ref{thmloc}]
Let $i:\cX\hookrightarrow \cX\times \CC$ denotes the inclusion defined by the graph of $f$. Identify $i_+\cM$ with $\cM[\partial _t]$. It is equipped with the good filtration
\begin{eqnarray*}
F_\ell\lefpar \cM[\partial _t]\rigpar &=&\sum_{j+k=\ell,\;k\geq 0}^{}F_j\cM\partial_{t}^{k}
\end{eqnarray*}
and $\gr_\ell^F\lefpar \cM[\partial_t]\rigpar $ is identified with $\oplus_{j\leq \ell}^{}\gr_j^F\!\cM$. The relative de~Rham complex
\begin{eqnarray*}
\DR_{\cX\times \CC/\CC}^{}(i_+\cM)&=&\lefpar \Omega_{\cX}^{\bbullet}\otimes_{\cO_\cX}^{}\cM[\partial _t],\nabla-\partial_t\cdot df\wedge\rigpar 
\end{eqnarray*}
is filtered by
\begin{eqnarray*}
F_p\DR_{\cX\times \CC/\CC}^{}(i_+\cM)&=&\lefpar \Omega_{\cX}^{\bbullet}\otimes_{\cO_\cX}^{}F_{p+\bbullet}^{}\cM[\partial _t],\nabla-\partial_t\cdot df\wedge\rigpar 
\end{eqnarray*}
so that
\begin{eqnarray*}
\gr_p^F\!\DR_{\cX\times \CC/\CC}^{}(i_+\cM)&=&\lefpar \Omega_{\cX}^{\bbullet}\otimes_{\cO_\cX}^{}\lefpar \oplus_{j\leq p+\bbullet}^{}\gr_j^F\!\cM\rigpar,\gr^F\!\nabla-df\wedge\rigpar.
\end{eqnarray*}
The right-hand term is also the $p$-th term of a filtration $G_\bbullet\lefpar \Omega_{\cX}^{\bbullet}\otimes_{\cO_\cX}^{}\gr^F\!\cM,\gr^F\!\nabla-df\wedge\rigpar$. The graded complex $G_p/G_{p-1}$ is the complex
\begin{eqnarray*}
\lefpar \Omega_{\cX}^{\bbullet}\otimes_{\cO_\cX}^{}\gr_{p+\bbullet}^F\!\cM,\gr^F\!\nabla\rigpar=\gr_p^F\!\DR\cM.
\end{eqnarray*}
If $p$ is large enough, this complex is acyclic in a neighbourhood of the compact fiber $f_{}^{-1}(0)\cap\supp\cM$ (see \eg. \cite{Laumon83,Malgrange85}). We conclude, taking inductive limits, that, for $p$ large enough and any $i$,
$$
\dim\bH^i\lefpar f_{}^{-1}(0),\lefpar \Omega_{\cX}^{\bbullet}\otimes_{\cO_\cX}^{}\gr^F\!\cM,\gr^F\!\nabla-df\wedge\rigpar\rigpar =\dim\bH^i\lefpar \cX\times \{0\},\gr_p^F\!\DR_{\cX\times \CC/\CC}^{}(i_+\cM)\rigpar .
$$

Let now $F_\bbullet\cH^if_+\cM$ be the good filtration defined as
$$
\mathop{\rm image}\nolimits \bR^if_*(F_p\DR_{\cX\times \CC/\CC}^{}(i_+\cM))\longrightarrow \bR^if_*(\DR_{\cX\times \CC/\CC}^{}(i_+\cM)).
$$
The strictness of the Hodge filtration on direct images (see \cite[th\ptbl 2.14]{MSaito87}) implies that
\begin{eqnarray*}
\dim\bH^i\lefpar \cX\times \{0\},\gr_p^F\!\DR_{\cX\times \CC/\CC}^{}(i_+\cM)\rigpar&=&\dim\gr_p^F\cH^if_+\cM
\end{eqnarray*}
for any $p$. Now, it follows from the local index theorem of Kashiwara (see \eg. \cite[p\ptbl67]{Malgrange91}) that, for $p$ large enough,
\begin{eqnarray*}
\dim\gr_p^F\cH^if_+\cM&=&\dim\phip_t\DR(\cH^if_+\cM)
\end{eqnarray*}
where $\phi_t$ denotes the vanishing cycle functor relative to $\id:\CC\rightarrow \CC$ and $\phip=\phi[-1]$. We conclude by using the same argument than at the end of \T\ref{num21}.
\end{proof}

\begin{proof}[Proof of theorem \protect\ref{thmglob}]
By theorem \ref{thmloc} and remark \ref{remphi}, it is enough to prove that the natural map
$$
\bH^i\lefpar X,\lefpar \Omega_{X}^{\bbullet}\otimes_{\cO_X}^{}\gr^F\!\cM,\gr^F\!\nabla-df\wedge\rigpar\rigpar\longrightarrow \bH^i\lefpar X_{}^{\rm an},\lefpar \Omega_{X}^{\bbullet}\otimes_{\cO_X}^{}\gr^F\!\cM,\gr^F\!\nabla-df\wedge\rigpar_{}^{\rm an}\rigpar\leqno{(*)}
$$
is an isomorphism for all $i$. The cohomology of the complex $\lefpar \Omega_{X}^{\bbullet}\otimes_{\cO_X}^{}\gr^F\!\cM,\gr^F\!\nabla-df\wedge\rigpar_{}^{\rm an}$ is supported by a finite number of fibers $f_{}^{-1}(c)$ (see \T\ref{numrem1}) and, by faithful flatness of $\cO_{X_{}^{\rm an}}^{}$ over $\cO_X$, the same holds for the complex $\lefpar \Omega_{X}^{\bbullet}\otimes_{\cO_X}^{}\gr^F\!\cM,\gr^F\!\nabla-df\wedge\rigpar$. These complexes can be viewed as having $\cO_{\wt X_{}^{\rm an}}^{}$- or $\cO_{\wt X}^{}$-coherent  cohomology (\cf. \T\ref{numnot} for the definition of $\wt X$), and GAGA implies that
$$
H^k\lefpar X,\cH^j\lefpar \Omega_{X}^{\bbullet}\otimes_{\cO_X}^{}\gr^F\!\cM,\gr^F\!\nabla-df\wedge\rigpar\rigpar \rightarrow H^k\lefpar X_{}^{\rm an},\cH^j\lefpar \Omega_{X}^{\bbullet}\otimes_{\cO_X}^{}\gr^F\!\cM,\gr^F\!\nabla-df\wedge\rigpar_{}^{\rm an}\rigpar
$$
is an isomorphism. The analytization morphism is compatible with the natural spectral sequences with $E_{2}^{jk}$-term described above. As $E_{2}^{jk}$ is finite dimensional, the spectral sequence degenerates at a finite rank and consequently $(*)$ is an isomorphism too.
\end{proof}

\providecommand{\bysame}{\leavevmode\hbox to3em{\hrulefill}\thinspace}

\end{document}